\documentclass{amsart}

\newcommand{\remark}{\noindent {\bf Remark. }}

\newcommand{\ws}{\hspace{4pt}}
\newtheorem{theorem}{Theorem}

\newtheorem{cor}{Corollary}
\newtheorem{lemma}{Lemma}
\newtheorem{e}{Example}
\newtheorem{defi}{Definition}
\usepackage[]{amssymb,amsopn}
\usepackage[]{amsmath}
\usepackage[]{eucal}
\usepackage[]{latexsym}
\usepackage{color}

\makeatletter
\@namedef{subjclassname@2020}{%
  \textup{2020} Mathematics Subject Classification}
\makeatother

\begin{document}

\title[]{Translation beyond Delsarte }
\author{\'A. P. Horv\'ath }

\thanks{Supported by the NKFIH-OTKA Grants K128922 and K132097.}

\begin{abstract} We introduce general translations as solutions to Cauchy or Dirichlet problems. This point of view allows us to handle for inctance the heat-diffusion semigroup as a translation. With the given examples Kolmogorov-Riesz characterization of compact sets in certain $L^p_\mu$ spaces are given. Pego-type characterizations are also derived. Finally for some examples the equivalence of the corresponding modulus of smoothness and K-functional is pointed out.

\end{abstract}
\maketitle

\section{Introduction}

In 1938 Delsarte introduced the notion of generalized translation, see \cite{d}. His starting point was the modification of Taylor's formula by the eigenfunctions of an ordinary differential operator: $D_xf(\lambda, \cdot)=\lambda f(\lambda,x)$. The Taylor series of $f=f(\lambda, z)$ around $z_0$ at $y$, with the notation $x=y-z_0$, can be expressed as
$$f(\lambda,x)=\sum_{k=0}^\infty \lambda^k \varphi_k(x),$$
where $D_x\varphi_0(x)=0$ and $D_x\varphi_k(x)=\varphi_{k-1}(x)$ ($\varphi_{-1}=0$). Moreover $\varphi_k(0)=0$ if $k>0$, and $\varphi_0(0)=1$. \\
He illustrated his idea with the following two examples:\\ $D_x(f)=f'$, $f(\lambda,x)=e^{\lambda x}$, $\varphi_k(x)=\frac{x^k}{k!}$ and\\  $D_x(f)=f''+\frac{2\alpha+1}{x}f'$, $f(\lambda,x)=j_\alpha(i\sqrt{\lambda} x)$, $\varphi_k(x)=\frac{\Gamma(\alpha+1)}{\Gamma(k+1)\Gamma(k+\alpha+1)}\left(\frac{x}{2}\right)^{2k}$;\\
and introduced the translation operator below:
$$T^t_x:=\sum_{k=0}^\infty\varphi_k(t)D_x^k.$$
The convergence of the series above implies that $T^t_xf(\lambda,\cdot)=f(\lambda,x)f(\lambda,t)$. Examination of the second example leads to the next chain of ideas. Denoting by $u(x,t):=T^t_xf$ the construction ensures that the translated functions are solutions to the next Cauchy problem.
$$D_tu=D_xu, \ws \ws u(x,0)= f(x), \ws \ws \frac{\partial}{\partial t}u(x,0)=0.$$
The initial study of the properties of generalized translation is due to Delsarte \cite{d} and Levitan see e.g. \cite{le}. Braaksma and Snoo in a series of papers dealt with the problem of introducing general translations via certain hyperbolic Cauchy problems, see e.g. \cite{b1}, \cite{bs}. Examination of product formulae with respect to classical orthogonal polynomials led to definition of Laguerre and Jacobi translations, see e.g.  \cite{g}, \cite{fjk1}, \cite{gm}.\\
General translation has widespread applications. Estimation of the operator norm leads to a maximum principle with respect to the hyperbolic problem in question. It gives certain convolution structures, see e.g. \cite{ag}, \cite{ka}. It also leads to estimation of $p$-Christoffel function, see e.g. \cite{adh}, \cite{abdh}. As in the standard case approximation theoretic problems, as smoothness, best approximation or approximation by Cesaro means can be examined by general translation, see e.g. \cite{gm1}, \cite{p}. By general translation compactness criteria can be derived too, see \cite{hor}.

Here we extend the notion of translation as solutions to Cauchy or Dirichlet problems in parabolic and elliptic cases as well. This approach allows to handle e.g. the heat-diffusion semigroup as a translation. Although several different examples are listed and investigated below, this one is highlighted, because unlike the old and the other new examples this one inherits the comfortable semigroup property of the standard translation, i.e. $T^{t_1}T^{t_2}=T^{t_1+t_2}$. Moreover parabolicity allows to draw conclusions in $d$-dimension while the variable of the translation is one-dimensional. We show that the general translations introduced below possess the properties of the standard one.

The paper is organized as follows. In the next section different types of examples of translation are given. In the third section the notion of "regularity with respect to compactness" is introduced and is applied to derive Kolmogorov-Riesz type characterizations in certain $L^p$-spaces. The rest of this section is devoted to show that all the listed translations are regular. In the fourth section Pego-type characterizations of compact sets of $L^2$ are derived by convolution method. The fifth section deals with the approximation theoretic aspect of translation that is with modulus of smoothness and K-functional.

\section{The translation operator}

We define the translation operator as follows.

\begin{defi}\label{D1} Let $\Omega=I\times J\subset \mathbb{R}^d\times \mathbb{R}_{+}$. Let $L$ be a linear partial differential operator of order at most two and take a function $f(x)\in C(I)$ or $f\in L^p_\mu(I)$, where $\mu$ is a Radon measure on $I$. Define translation as the next linear operator: the translation of $f$
\begin{equation} T_x^tf:=T^tf(x)=u(x,t),\end{equation}
where $u(x,t)$ is the solution to
\begin{equation}\label{ade} Lu=0, \ws (x,t)\in \Omega,\end{equation}
\begin{equation}\label{kef}  u(x,0)=f(x), \ws x\in I,\end{equation}
where the last equality is meant in $\sup$-norm or in $p$-norm.
\end{defi}

\remark

All of our examples can be given by appropriate integral transformations as well. The corresponding kernel functions are denoted by $W_t(x,y)$, $K_t(x,y)$, etc. according to the standard notation. Here $t>0$, and $x,y\in \mathbb{R}^d$, $d\ge 1$. That is
\begin{equation}\label{trdef} T^tf(x):=\int_IW_t(x,y)f(y)dy.\end{equation}
It also makes the solution well-defined.

\medskip

Below different types of examples are introduced. Of course, by the same chain of ideas several further examples can be constructed. Here are the ones we study subsequently.

\medskip

\subsection{Heat-diffusion semigroup associated with Hermite functions} This semigroup, from different point of views, is investigated by several authors, see e.g. \cite{st}, \cite{at}, \cite{gs} and the references therein.

In $\mathbb{R}^d$ the eigenfunctions of the $d$-dimensional harmonic oscillator (Hermite operator),
\begin{equation}\label{ho} D_x=\Delta_x-|x|^2,\end{equation}
are the $d$-dimensional Hermite functions,
$$\tilde{h}_\nu=\tilde{h}_{n_1}\dots \tilde{h}_{n_d},$$
where $\nu=(n_1, \dots ,n_d)$ $n_i\in\mathbb{N}$ and
$$\tilde{h}_{k}(x)=h_k(x)e^{-\frac{x^2}{2}}=\frac{1}{\sqrt{2^k k! \sqrt{\pi}}}H_k(x)e^{-\frac{x^2}{2}}.$$
($H_k$ are the Hermite polynomials cf. \cite{sz}.)
The associated heat-diffusion semigroup is given by its kernel function defined on $\mathbb{R}^d\times \mathbb{R}^d \times \mathbb{R}_+$,
\begin{equation}\label{hsk}W_t(x,y)=\sum_{n=0}^\infty e^{-(2n+d)t}\sum_{|\nu|=n}\tilde{h}_\nu(x)\tilde{h}_\nu(y).\end{equation}
Thus for and appropriate $f$, denoting by $\tilde{f}(y):=f(y)e^{-\frac{y^2}{2}}$,
\begin{e}\label{her}$$T^tf(x):=\int_{\mathbb{R}^d}W_t(x,y)\tilde{f}(y)dy.$$\end{e}
Denoting by $u(x,t):=T^tf(x)$, if $\tilde{f}\in L^p(\mathbb{R}^d)$, $1\le p \le\infty$, then
\begin{equation}\label{hde}\left(\frac{\partial}{\partial t}-D_x\right)u(x,t)=0\end{equation}
and
\begin{equation}\lim_{t\to 0+}\|T^tf-\tilde{f}\|_p=0,\end{equation}
cf. \cite[Proposition 2.5 and Theorem 2.6]{st}.

\medskip

\begin{center} -- $\cdot$ -- \end{center}

\medskip

\subsection{Further Poisson integrals - elliptic and parabolic equations} Following the previous track of thoughts translation of a function in $f\in L^p(\mathbb{R})$ can be defined by  Poisson integrals. First we take the simplest elliptic and parabolic equations on the upper half-plane and define translations by standard convolution.

In the next example $T_x^tf$ is a harmonic function on the (open) upper half-plane with limit $f$ on the real line and with limit zero at infinity cf. e.g. \cite{kh}, that is the solution to this Dirichlet problem
$$\Delta u(x,t)=0, \ws \ws(x,t) \in  \mathbb{R}\times \mathbb{R}_+, \ws \ws u(x,0)=f(x).$$
It is given by the Poisson integral below.

\begin{e}\label{th} Let $f\in L^p(\mathbb{R})$, $1\le p <\infty$ or $f\in C_b(\mathbb{R})$.
\begin{equation}\label{htr} T^tf(x):=\frac{1}{\pi}\int_\mathbb{R}f(\xi)\frac{t}{(\xi-x)^2+t^2}dt=f*K_{t,e}(x),\end{equation}\end{e}
where $C_b(\mathbb{R})$ stands for continuous and bounded functions on the real line.

\medskip

\begin{center} -- $\cdot$ -- \end{center}

\medskip

Considering the Cauchy problem
$$u_t(x,t)=b^2\Delta u(x,t) \ws \ws (x,t)\in\mathbb{R}\times (0,\infty), \ws u(x,0)=f(x),$$
the next translation can be defined by the corresponding Poisson integral.
\begin{e}\label{tp} Let $f\in L^p(\mathbb{R})$ ($1\le p \le\infty$).
\begin{equation}\label{ptr} T^tf(x):=\int_\mathbb{R}f(\xi)\frac{e^{-\frac{(x-\xi)^2}{4b^2t}}}{2b\sqrt{\pi}\sqrt{t}}dt=f*K_{t,p}(x).\end{equation}\end{e}

\medskip

\begin{center} -- $\cdot$ -- \end{center}

\medskip

We continue with an elliptic example again, where the Poisson integral is not of convolution type.\\
Let $P_n^{(\lambda)}(\cos \vartheta)$ be the $n$th ultraspherical polynomial ($\lambda>-\frac{1}{2}$, $\vartheta \in [0,\pi]$) orthogonal with respect to $d\mu(\vartheta)=\sin^{2\lambda}\vartheta d\vartheta$, see e.g. \cite{sz}. An $f \in L^p_\mu$ with some $1 \le p \le \infty([0,\pi])$ has an expansion $f\sim \sum_{k=0}^\infty a_k P_k^{(\lambda)}(\cos \vartheta)$ with $a_k=\gamma_k\int_0^\pi f(\vartheta)P_k^{(\lambda)}(\cos \vartheta)d\mu(\vartheta)$. In \cite{ms} the authors examined the next Poisson integral of $f$,
$$u(x,t)=f(r,\vartheta)=\sum_{k=0}^\infty a_k r^kP_k^{(\lambda)}(\cos \vartheta)=\int_0^\pi P(r,\vartheta,\varphi) f(\varphi)d\mu(\varphi),$$
where $0<r<1$, $x=r\cos\vartheta$, $t=r\sin\vartheta$. The Poisson kernel is
$$P(r,\vartheta,\varphi)=\sum_{k=0}^\infty \gamma_k r^k P_k^{(\lambda)}(\cos \vartheta)P_k^{(\lambda)}(\cos \varphi).$$
Then $u(x,t)$ satisfies the differential equation
$$u_{xx}+u_{tt}+\frac{2\lambda}{t}u_t=0.$$

Since $\|f(r,\vartheta)-f(\vartheta)\|_{\mu,p}\to 0$ as $r \to 1$ if $1\le p<\infty$ and also for $p=\infty$ if $f$ is continuous, in view of Definition \ref{D1} we have the next example. As above, the closed form of the Poisson kernel implies the definition below.

\begin{e}\label{tm} Let $D=[0,\pi]\times [0,1)$,  $f \in L^p_\mu([0,\pi])$.
$$T^tf(x)=u(x,t)=T^rf(\vartheta)=f(r,\vartheta)=\int_0^\pi P(r,\vartheta,\varphi) f(\varphi)d\mu(\varphi),$$
where
\begin{equation}\label{p}P(r,\vartheta,\varphi)=\frac{\lambda(1-r^2)}{\pi}\int_0^\pi\frac{\sin^{2\lambda-1}\xi}{\left(1-2r(\cos\vartheta\cos\varphi+\sin\vartheta\sin\varphi\cos \xi)+r^2\right)^{\lambda+1}}d\xi.\end{equation}
\end{e}

\medskip

\begin{center} -- $\cdot$ -- \end{center}

\medskip

Considering the next differential operator
$$D_x=L^{\alpha}:=\frac{1}{2}\left(-\frac{d^2}{dx^2}+x^2+\frac{\alpha^2-\frac{1}{4}}{x^2}\right)$$
with eigenfunctions
$$L^{\alpha}\varphi_n^{(\alpha)}=(2n+\alpha+1)\varphi_n^{(\alpha)},$$
where
\begin{equation}\label{fi}\varphi_n^{(\alpha)}(x)=\sqrt{\frac{\Gamma(n+1)}{\Gamma(n+\alpha+1)}}e^{-\frac{x^2}{2}}L_n^{(\alpha)}(x^2)\sqrt{2x}\end{equation}
($L_n^{(\alpha)}$ is the standard Laguerre polynomial, see e.g. \cite{sz}), we define a translation by the solution to the parabolic Cauchy problem
$$u_t+L^{\alpha}u=u_t-\frac{1}{2}u_{xx}+\frac{1}{2}\left(x^2+\frac{\alpha^2-\frac{1}{4}}{x^2}\right)u=0, \ws \ws u(x,0)=g(x),$$
$x,t\in (0,\infty)$,  $\alpha\ge -\frac{1}{2}$, cf. \cite{ch}.

\begin{e}\label{tpl}
$$T^tg(x):=\int_0^\infty W_t^\alpha(x,y)g(y)dy,$$
where
$$W_t^\alpha(x,y)=2\sqrt{xy}\frac{e^{-t}}{1-e^{-2t}}I_\alpha\left(\frac{2xye^{-t}}{1-e^{-2t}}\right)e^{-\frac{1}{2}(x^2+y^2)\frac{1+e^{-2t}}{1-e^{-2t}}},$$
\end{e}
where $I_\alpha(z)=\sum_{k=0}^\infty \frac{\left(\frac{z}{2}\right)^{2k+\alpha}}{\Gamma(k+1)\Gamma(k+\alpha+1)}$ is the modified Bessel function.

\medskip

\begin{center} -- $\cdot$ -- \end{center}

\medskip

\subsection{Fourier method - hyperbolic equation} The first - well known - examples are given by hyperbolic equations generated by Sturm-Liouville type operators, see e.g. \cite{d}, \cite{le}. Consider the operator
$$D_x=\frac{\partial^2}{\partial x^2}+q(x)\frac{\partial}{\partial x}-r(x).$$
If
$$D_x\varphi(x)=a\varphi(x),$$
denoting by $u(x,t)=\varphi(x)\varphi(t)$, $u(x,t)$ fulfils the equation below.
\begin{equation}\label{Lh} u_{xx}-u_{tt}+q(x)u_x-q(t)u_t-(r(x)-r(t))u=0.\end{equation}

In Laguerre and Bessel cases $q(x)=\frac{2\alpha+1}{x}$, where $\alpha\ge -\frac{1}{2}$ and $r(x)=x^2$ or $r(x)=0$, respectively. The eigenfunctions of $D_x$ are the Laguerre functions $\mathcal{L}_n^{(\alpha)}(x)=\frac{n!\Gamma(\alpha+1)}{\Gamma(n+\alpha+1)}L_n^{(\alpha)}(x^2)e^{-\frac{x^2}{2}}$, and the Bessel functions $j_\alpha(\lambda x)$, where $L_n^{(\alpha)}(x)$ and $j_\alpha(\lambda x)$ stand for the classical Laguerre polynomials and the entire Bessel functions, respectively. These examples are studied by several authors, see e.g. \cite{ag}, \cite{adh}, \cite{gm}, \cite{p}, \cite{abdh}.

In Jacobi case with  $\alpha \ge\beta \ge -\frac{1}{2}$, $\alpha > -\frac{1}{2}$ by the same argument,
 $q(x)=\frac{(\beta+1)\cos x-\beta}{\sin x}$ and $r(x)=\frac{\alpha(\alpha-2\beta+1)\cos x}{4(1-\cos x)}$. The eigenfunctions are the Jacobi functions, $\mathcal{P}^{(\alpha,\beta)}_n(x)=c_nP^{(\alpha,\beta)}_n(\cos x)\sin^\alpha\frac{x}{2}$, where $P_n^{(\alpha, \beta)}(x)$-s are the classical Jacobi polynomials. For Jacobi translation see e.g. \cite{g}, \cite{ad}.

 In the listed cases the eigenfunctions of the Sturm-Liouville operators $\{u_k\}_{k=0}^\infty$ form Schauder bases in the corresponding (weighted) spaces, that is if the initial condition is given by $f(x)\sim \sum_{k=0}^\infty a_ku_k(x)$, then the solution to the Cauchy problem can be given as $T^tf(x)=u(x,t)\sim \sum_{k=0}^\infty a_ku_k(x)u_k(t).$

Thus translation generated by a hyperbolic equation is as follows.
$$T^tf(x)=u(x,t)$$
if
\begin{equation}\label{BLJ} L_hu=0, \ws\ws (x,t)\in I\times I; \ws \ws u(x,0)=f(x), \ws \ws u_t(x,0)=0,\ws x\in I.\end{equation}
Here $L_h$ is given by \eqref{Lh}, where $I=\mathbb{R}_+$ in Laguerre and Bessel cases and in Jacobi case $I=(0,\pi)$.

It is pointed out in a series of papers that Bessel, Laguerre and Jacobi translations are bounded operators, see the cited papers above and the references therein. The method described above implies the symmetry of the translation derived by a Sturm-Liouville equation, see \cite{b1}, \cite{bs}.\\

Next we give an example by Jacobi functions. Let $\alpha \ge\beta \ge -\frac{1}{2}$, $\alpha > -\frac{1}{2}$ again and let us denote by $\varrho=\alpha+\beta+1$ and now let $\tau \in \mathbb{R}$.
Let
\begin{equation}\label{w} w(x):=2^{2\varrho}\sinh^{2\alpha+1}x\cosh^{2\beta+1}x.\end{equation}
With $q(x)=\frac{w'}{w}(x)$, $\varphi_\tau^{(\alpha,\beta)}(x)$ fulfils the differential equation:
$$y''+qy'+(\varrho^2+\tau^2)y=0,$$
where
$$\varphi_\tau^{(\alpha,\beta)}(x):=\varphi_\tau(x)=_2F_1\left(\frac{\varrho+i\tau}{2},\frac{\varrho-i\tau}{2};\alpha+1;-\sinh^2x\right)$$
are the Jacobi functions. Thus the translation is defined as it follows, cf. e.g. \cite{p1}, \cite{p2}, \cite{git} and the references therein.
\begin{e}\label{jfv} Let $f$ be a suitable even function on $\mathbb{R}$. $T^tf(x)=u(x,t)$ if
$$u_{xx}-u_{tt}+q(x)u_x-q(t)u_t=0, \ws \ws u(x,0)=f(x), \ws \ws u_t(x,0)=0,\ws x\in \mathbb{R},$$
where
$$q(x)=(2\alpha+1)\coth x+(2\beta+1)\tanh x.$$
\end{e}

\medskip

\section{Regularity with respect to compactness}

In this section we derive Kolmogorov-Riesz-type compactness criteria in certain $L^p$-spaces. For the original theorem see e.g. \cite{a}. This theorem has several extensions to different function spaces with standard translation, see e.g. \cite{au}, \cite{dfg}, \cite{gma}, \cite{gr}, \cite{r}, \cite{hohm} and with Bessel and Laguerre translations, see \cite{hor}. In all but one of our examples the translation is not symmetric. This implies that one of the criteria of compactness has to be splitted into two different assumptions, see Definition \ref{eqcdef} below.

\medskip

We start with some notation. Let $\mu$ and $\nu$ be Radon measures on $I$ and $J$, respectively, and $1\le p<\infty$. $\|f\|_{\mu(x),p}=\left(\int_I|f|^pd\mu\right)^{\frac{1}{p}}$. $T^0f(x)=f(x)$.

We define the "norm" of the translation operator $M_T$ as
\begin{equation}\label{Tnorm} M_T:=M_{T,p}=\sup_{t}\sup_{\|f\|_{\mu,p}\le 1}\|T^tf(x)\|_{\mu(x),p}, \ws \ws 1\le p \le \infty.\end{equation}

\medskip

We introduce the next notation: Let $a>0$, $M_0>0$ fixed, $B_a:=\{\|x\|\le a\}\cap I$, or $B_a:=(-a,a)\cap J$; $A:=\int_{B_a}1d\nu(t)$.
\begin{equation}M_{a,R}f(x):=\left\{\begin{array}{ll}\frac{1}{A}\int_{B_a}T^tf(x)d\nu(t) \ws |x|<R,\\ 0 \ws |x|\ge R.\end{array}\right.\end{equation}

\medskip

\begin{defi}\label{regd} A translation, $T^t_x$ is regular with respect to $p$, $\mu$, $\nu$ if it fulfils the properties below. \\
{\bf \rm{(T1)}} There is a dense set $\mathcal{E}_1\subset L^p_\mu(I)$ such that for all $g\in \mathcal{E}_1$, $0<t<M_0$ and $\varepsilon>0$ there is a $\delta=\delta(\varepsilon, M_0, g)>0$ such that for all $0\le |h| \le \delta$
\begin{equation}\label{leqf}\left(\int_I \left|\left(T^{t+h} g(x)-T^{t} g(x)\right)\right|^pd\mu(x)\right)^{\frac{1}{p}} <\varepsilon. \end{equation}
{\bf \rm{(T2)}} There is a dense set $\mathcal{E}_2\subset L^p_\mu(I)$ such that for all $g\in \mathcal{E}_2$, and positive numbers $\varepsilon$, $a$, $R$ there is a $\delta=\delta(g,\varepsilon,a,R)>0$ such that if $0\le |h| \le \delta$, then
\begin{equation}\label{leqM} \left|M_{a,R}g(x+h)-M_{a,R}g(x)\right| <\varepsilon,  \ws \ws x, x+h \in B_R. \end{equation}
{\bf \rm{(T3)}} There is a finite constant, $c(a,R)$ so that for all $f \subset L^p_\mu(I)$,
\begin{equation}\label{fy} |M_{a,R}f(x)| \le c(a,R)\|f\|_{\mu,p}, \ws \ws \forall \ws x\in I.\end{equation}
\end{defi}

\medskip

\remark

\noindent (R1) If $T_x^y=T_y^x$ and $\mu=\nu$, then \eqref{leqf} implies  \eqref{leqM} and if $M_T<\infty$ it also implies \eqref{fy}. Indeed
$$\left|M_ag(x+h)-M_ag(x)\right|\le \frac{1}{A} \left(\int_{B_a}1d\mu(y)\right)^{\frac{1}{p'}}\left(\int_{B_a}|T^{x+h}g(y)-T^{x} g(y)|^pd\mu(y)\right)^{\frac{1}{p}},$$
which ensures \eqref{leqM}. \eqref{fy} can be obtained in the same way.

\medskip

\begin{defi} A set $K\subset L^p_\mu(I)$ is equivanishing, or we say it fulfils property ${\bf P_a}$ if
 for all $\varepsilon >0$ there is an $R>0$ such that  for all $f\in K$
$$ \left(\int_{I\setminus B_R} |f(x)|^pd\mu(x)\right)^{\frac{1}{p}}<\varepsilon. $$
\end{defi}

\medskip

\begin{defi}\label{eqcdef} A set $K\subset L^p_\mu(I)$ is equicontinuous in mean if it fulfils the next properties\\
${\bf P_{b_1}}$: For all $\varepsilon$ and $M_0$ positive numbers there is a $\delta>0$ (independent of $f$) such that for all $t\in B_{M_0}$, $0\le |h| \le \delta$ and $f\in K$ \eqref{leqf} is satisfied.\\
${\bf P_{b_2}}$: For all $\varepsilon>0$ there is a $\delta=\delta(a,R)>0$ (independent of $f$) such that for all $0\le |h| \le \delta$ and $f\in K$ \eqref{leqM} is satisfied.
\end{defi}

\medskip

\begin{theorem}\label{T1} With the notation above, let $T=T^t_x$ be a bounded translation in \eqref{Tnorm}-sense and suppose that it is regular with respect to $p$, $\mu$, $\nu$. Let $K\subset L^p_\mu(I)$ be a bounded set. Then $K$ is precompact if and only if it is equivanishing and equicontinuous in mean. \end{theorem}

\proof
First we prove that if $K$ is bounded, equivanishing and equicontinuous in mean and the translation has the property {\bf T3}, then $K$ is precompact.\\
For an arbitrary $\varepsilon>0$ let $R$ be chosen according to ${\bf P_a}$. Let $a>0$ be fixed and will be chosen later, and define
$$F_{a,R}:=\{M_{a, R}f(x): f\in K, \ws x\in B_R\}.$$
Then {\bf T3} and the boundedness of $K$ imply that $F_{a,R}$ is (uniformly) bounded. According to ${\bf P_{b_2}}$ $F_{a,R}$ is equicontinuous, thus it is precompact. Let $M_af_1, \dots , M_af_n$ be an $\varepsilon_1$-net in $F_{a,R}$, where $\varepsilon_1=\frac{\varepsilon}{(\mu(B_R))^{\frac{1}{p}}}$. We show that $f_1, \dots ,f_n$ is a $5\varepsilon$-net in $K$.
Recalling the definition of $R$ and $A$
$$\| M_{a,R}f-f\|_{\mu,p}\le \left\| M_{a,R}f-f|_{B_R}\right\|_{\mu,p} +\left(\int_{I\setminus B_R}\left|f\right|^pd\mu\right)^{\frac{1}{p}}$$
$$\le \left(\int_{B_R}\left|\frac{1}{A}\int_{B_a}\left(T^{t}f(x)-f(x)\right)d\nu(t)\right|^pd\mu(x)\right)^{\frac{1}{p}}+ \varepsilon $$
$$\le \frac{1}{A}\int_{B_a}\left(\int_I\left|T^{t}f(x)-f(x)\right|^pd\mu(x)\right)^{\frac{1}{p}}d\nu(t)+ \varepsilon$$
$$\le \sup_{0\le |t| \le a}\left\|T^{t}f(x)-f(x)\right\|_{\mu,p}+ \varepsilon.$$
Considering ${\bf P_{b_1}}$, if $a$ is small enough, $\| M_{a,R}f-f\|_{\mu,p}\le 2\varepsilon$ for all $f\in K$. Let $f\in K$ be arbitrary. Selecting $f_i$ such that $|M_{a,R}f(x)-M_{a,R}f_i(x)|\le\varepsilon_1$,
$$\|M_{a,R}f-M_{a,R}f_i\|_{\mu,p}\le \varepsilon_1(\mu(B_R))^{\frac{1}{p}}=\varepsilon.$$
Thus by triangle inequality
$$\|f-f_i\|_{\mu,p}\le 5\varepsilon.$$

On the other hand suppose that $K$ is precompact. Let $\varepsilon>0$ be arbitrary and $\varphi_1, \dots ,\varphi_n$ be an $\frac{\varepsilon}{2}$-net in $K$. As $C_0$, that is the set of compactly supported continuous functions, is dense in $L^p_\mu(I)$, there are $g_1, \dots , g_n \in C_0$ so that for all $f \in K$ there is an $i$ such that $\|g_i-f\|_{\mu,p}<\varepsilon$. Then $R$ is appropriate with respect to ${\bf P_a}$ if $B_R$ contains the supports of $g_i$, $i=1,\dots n$.\\
Since $T$ is bounded, selecting $g_1, \dots , g_n \in \mathcal{E}_1$ as above, triangle inequality together with {\bf T1} ensures ${\bf P_{b_1}}$.\\
Similarly, selecting $g_1, \dots , g_n \in \mathcal{E}_2$ as above, {\bf T2}, {\bf T3} and the triangle inequality imply ${\bf P_{b_2}}$.

\medskip

\remark  (R2) For compactness it is enough to assume that the initial condition, \eqref{kef} is fulfilled uniformly. The importance of the extension above will be shown in the next section.

\medskip

In the rest of this section we show that all the translations listed above are regular in sense of Definition \ref{regd}. At first we make some observations which will be useful in parabolic cases below and in the last section.

\medskip

\remark (R3) For the next observation let us write the differential equations of the listed examples in the next form:
\begin{equation}\label{des}D_{(1),t}u(x,t)=D_{(2),x}u(x,t), \ws \ws u(x,0)=f(x). \end{equation}
$T^tf(x)=u(x,t)$ can be defined by an integral transformation, that is
$$T^tf(x)=\int_If(\xi)K(x,t,\xi)w(\xi)d\xi.$$
Thus the kernel function fulfils the differential equation above. As the kernel function is symmetric at least in $x$ and $\xi$, i.e.
$$K(x,t,\xi)=K(\xi,t,x),$$
we have
$$D_{(1),t}K(x,t,\xi)=D_{(2),x}K(x,t,\xi)=D_{(2),\xi}K(x,t,\xi).$$
These observations imply the next lemma.

\medskip

\begin{lemma} Let $g$ be a smooth function on $I$ and assume that its derivatives disappears quickly enough at the boundary of $I$. Then in the listed cases
\begin{equation}D_{(2)}T^tg=T^tD_{(2)}g.\end{equation}
\end{lemma}

\proof Recall, that in our one dimensional examples
$$D_{(2)}h=h''+\frac{w'}{w}h'+rh.$$
Integrating by parts and considering the boundary condition we have
$$D_{(2),x}T^tg(x)=\int_Ig(\xi)D_{(2),x}K(x,t,\xi)w(\xi)d\xi=\int_Ig(\xi)D_{(2),\xi}K(x,t,\xi)w(\xi)d\xi$$ $$=\int_I-\frac{\partial}{\partial \xi}K(x,t,\xi)(gw)'(\xi)+(gw)(\xi)\left(\frac{w'(\xi)}{w(\xi)}\frac{\partial}{\partial \xi}K(x,t,\xi)+r(\xi)\right)w(\xi)d\xi$$ $$=\int_IK(x,t,\xi)D_{(2),\xi}g(\xi)w(\xi)d\xi=T^tD_{(2)}g(x).$$
If $w\equiv 1$, it is more direct. Let us see for instance the heat diffusion semigroup. Recalling that $I=\mathbb{R}^d$
$$\int_Ig(\xi)\left(\Delta_\xi W_t(x,\xi)-|\xi|^2W_t(x,\xi)\right)d\xi=\int_I\left(\Delta_\xi g(\xi)-g(\xi)|\xi|^2\right)W_t(x,\xi)d\xi.$$

\medskip

\begin{cor} Let $D_{(1)}=\frac{\partial}{\partial t}$ and let $g$ be as above. Then
\begin{equation}\label{delta} \|T^{t+h}g-T^tg\|_{p,\mu}\le h M_T\|D_{(2)}g\|_{p,\mu}.\end{equation}
\end{cor}

\proof
$$T^{t+h}g(x)-T^tg(x)=h\frac{\partial}{\partial t}T^{t+\eta}g(x)=hD_{(2),x}T^{t+\eta}g(x)=hT^{t+\eta}D_{(2)}g(x).$$
Taking into consideration \eqref{Tnorm} we get the estimation.

\medskip

\begin{lemma} The translation given in Example \ref{her} is regular and $M_T=1$.
\end{lemma}

\proof The boundedness of the operator by a constant $C$ is proved in \cite[Theorem 2.6]{st}. It can be shown similarly that $M_T=1$. Indeed, according to Mehler's formula and \cite[Remark 2.5]{gs}
$$0<W_t(x,y)=(2\pi \sinh 2t)^{-\frac{d}{2}}e^{-\frac{|x-y|^2}{2}\coth 2t-\langle x, y\rangle \tanh t}$$ \begin{equation}\label{magbecs}\le (2\pi t)^{-\frac{d}{2}}e^{-\frac{|x-y|^2}{4t}}=K_t(x-y).\end{equation}
Thus
$$|T^tf|\le |\tilde{f}|*K_t$$
and since $\|K_t\|_1=1$
$$\|T^tf\|_p\le \| |\tilde{f}|*K_t\|_p\le \|\tilde{f}\|_p,$$
$1\le p \le \infty$.

Let $\mathcal{E}_1= C^2_0\subset L^p(\mathbb{R}^d)$. In view of \eqref{delta} {\bf \rm{(T1)}} is fulfilled.

We prove property {\bf \rm{(T2)}} by similar arguments. Considering the uniform convergence of \eqref{hsk} for a $t>0$, the recurrence and derivation formulae for Hermite functions imply that
$$\frac{\partial}{\partial x_j}W_t(x,y)=(y_je^{-2t}-x_j)W_t(x,y)-e^{-2t}\frac{\partial}{\partial y_j}W_t(x,y).$$
Thus
$$W_t(x+\chi,y)-W_t(x,y)=\langle e^{-2t}y-(x+\xi),\chi \rangle W_t(x+\xi,y)-e^{-2t}\langle \nabla_yW_t(x+\xi,y), \chi\rangle.$$
Let $f\in C^1_0(\mathbb{R}^d)$, $\mathrm{supp}f\subset B_r$, $x\in B_R$.
$$\left|\frac{1}{a}\int_0^{a}\int_{\mathbb{R}^d}(W_t(x+\chi,y)-W_t(x,y))\tilde{f}(y)dydt\right|$$ $$\le \frac{1}{a}\int_0^{a}\left(2|\chi|(r+R)\|T^t|f|\|_p+\int_{\mathbb{R}^d}W_t(x+\chi,y)\langle \chi, \nabla_y \tilde{f}(y)\rangle dy\right)dt$$ $$\le |\chi|\left(2(r+R)\|\tilde{f}\|_p+\|\nabla_y \tilde{f}\|_p\right).$$

To prove {\bf \rm{(T3)}} let us recalling \eqref{magbecs}. Then we have
$$\left|\frac{1}{a}\int_0^{a}T^tfdt\right|\le \frac{1}{a}\int_0^{a}\int_{\mathbb{R}^d}|\tilde{f}(y)|K_t(x-y)dydt\le \|\tilde{f}\|_p.$$

\begin{lemma} The translation given in Example \ref{jfv} is regular and $M_T=1$.
\end{lemma}

\proof For boundedness see \cite[Lemma 5.2]{fjk1}.\\
As it is mentioned above, in \cite{bs} it is pointed out that a translation generated by a hyperbolic equation is necessarily symmetric. It can be expressed by an integral transform whose kernel is symmetric in its three variables, see \cite{fjk1}. Thus according to remark (R1) it is enough to prove property {\bf (T1)}. Since $f$ is even, we can take $I=\mathbb{R}_+$. Let $\mathcal{D}=C^1_0(\mathbb{R}_+)\subset L^p_\mu (\mathbb{R}_+)$ the dense set in question, where
\begin{equation}\label{mu}d\mu(x)=\frac{1}{\sqrt{2\pi}}w(x)dx,\end{equation}
cf. \eqref{w}.
Let $f\in C^1_0(\mathbb{R}_+)$ and let $f(z)=g(\cosh 2z)$. In view of \cite[(5.1)]{fjk1} and \cite[(4.16)]{fjk1}
$$T^tf(x)=\int_{\mathbb{R}_+}f(z)K(x,t,z)d\mu(z)=\int_0^1\int_0^\pi g(x,t,r,\psi)dm(r,\psi),$$
where
$$(x,t,r,\psi)$$ $$=2(\cosh^2 x\cosh^2 t+r^2\sinh^2 x \sinh^2 t +2r\cos\psi\sinh x\sinh t \cosh x\cosh t)-1,$$
and
$$dm(r,\psi)=\frac{2\Gamma(\alpha+1)}{\sqrt{\pi}\Gamma(\alpha-\beta)\Gamma(\beta+\frac{1}{2})}(1-r^2)^{\alpha-\beta-1}r^{2\beta+1}\sin^{2\beta}\psi drd\psi.$$
It can be readily seen that
$$2\cosh^2(x-t)-1\le(x,t,r,\psi)\le 2\cosh^2(x+t)-1.$$
Since $g\in C^1_0(\mathbb{R}_+)$ too, suppose that $\mathrm{supp}g(z)\subset B_R$. Recalling that $0<t<M_0$ $\mathrm{supp}g_{t,r,\psi}(x)\subset B_{ C(M_0,R)}$. Considering that $\frac{\partial}{\partial y}g(x,t,r,\psi)$ is bounded and $m(r,\psi)$ is a bounded measure one can conclude that {\bf (T1)} fulfils with $\delta= c \varepsilon$.

\medskip

\medskip

\begin{lemma} The translations given by Examples \ref{th} and \ref{tp} are regular and $M_T=1$ in both cases.
\end{lemma}

\proof
The proof is based on the standard convolution structure.\\
First let us observe that the kernel functions are positive and $\int_\mathbb{R}K_{e,p}(x,t)dx=1$ for all $t>0$, which gives the operator norm.\\
In both cases let $\mathcal{E}_1=\mathcal{E}_2=C_0(\mathbb{R})$.
In harmonic case by the standard substitution $\frac{\xi-x}{t+h}=\tan \varphi$, where $h>0$ and $h=0$ in the different integrals, we have
$$\left\|T^{t+h} g(x)-T^{t} g(x)\right\|_{p}$$ $$=\left(\int_\mathbb{R}\left|\frac{1}{\pi}\int_{-\frac{\pi}{2}}^{\frac{\pi}{2}}(g(x+(t+h)\tan\varphi)-g(x+t\tan\varphi))d\varphi\right|^pdx\right)^{\frac{1}{p}}$$ $$\leq \left(\int_\mathbb{R}\left|\frac{1}{\pi}\int_{-\frac{\pi}{2}}^{-\alpha}(\cdot))d\varphi\right|^pdx\right)^{\frac{1}{p}}+\left(\int_\mathbb{R}\left|\frac{1}{\pi}\int_{-\alpha}^{\alpha}(\cdot))d\varphi\right|^pdx\right)^{\frac{1}{p}}+\left(\int_\mathbb{R}\left|\frac{1}{\pi}\int_{\alpha}^{\frac{\pi}{2}}(\cdot))d\varphi\right|^pdx\right)^{\frac{1}{p}}$$ $$=I+II+II.$$
$$I, \ws III \le 2 \|g\|_p\frac{\frac{\pi}{2}-\alpha}{\pi},$$
which is small if $\alpha$ is close enough to $\frac{\pi}{2}$.\\
Since $g$ is uniformly continuous for fixed $\alpha$ and $\varepsilon$ we can choose $h$, $|h|<1$ such that $|h\tan\varphi|$ is small enough to be the integrand in $II$ small, less than $\varepsilon$, say.\\
Then recalling that $g$ is compactly supported and $|t| \le M_0$, there is an $R$ such that the support of $g(x+(t+h)\tan\varphi)-g(x+t\tan\varphi)$ is in $[-R,R]$. Thus $II\le\varepsilon (2R)^{\frac{1}{p}}.$

Similarly in parabolic case let $u=\frac{\xi-x}{2b\sqrt{t+h}}$. Then
$$\left\|T^{t+h} g(x)-T^{t} g(x)\right\|_{p}$$ $$=\left(\int_\mathbb{R}\left|\frac{1}{\sqrt{\pi}}\int_{\mathbb{R}}(g(x+2b\sqrt{t+h}u)-g(x+2b\sqrt{t}u))e^{-u^2}du\right|^pdx\right)^{\frac{1}{p}}$$ $$\le \frac{1}{\sqrt{\pi}}\int_{|u|>R}e^{-u^2}\left(\int_\mathbb{R}|g(x+2b\sqrt{t+h}u)-g(x+2b\sqrt{t}u)|^pdx\right)^{\frac{1}{p}}du+\frac{1}{\sqrt{\pi}}\int_{|u|\le R}(\cdot)$$ $$=I+II.$$
Recalling that $g$ is compactly supported, the $p$-norm inside is bounded, thus if $R$ is large enough, $I$ is small. The uniform continuity of $g$ implies that $II$ is small if $h$ is small enough.

 Turning to the proof of \eqref{leqM} and \eqref{fy} we define $\nu(t)$. If $1<p<\infty$ and in parabolic case also for $p=1$ let $d\nu(t)=dt$, and if $p=1$, in the elliptic case let $d\nu(t)=\sqrt{t}dt$. Then the same replacements imply the results. Indeed, let us see the elliptic case first.
$$|M_ag(x+h))-M_ag(x)|$$ $$\le \frac{1}{A\pi}\int_{B_a} \left(\int_{-\frac{\pi}{2}}^{-\alpha}(\cdot)+\int_{-\alpha}^{\alpha}|g(x+(t+h)\tan\varphi)-g(x+t\tan\varphi)|d\varphi + \int_{\alpha}^{\frac{\pi}{2}}(\cdot)\right)d\nu(t).$$
In the first and third integrals the boundedness of the integrand, in the second one the uniform continuity of $g$ ensures the required estimates.\\
To prove \eqref{fy} let $1<p<\infty$.
$$|M_af(x)|\le \frac{1}{A}\int_{B_a}\|f\|_p\frac{1}{\pi t}\left(\int_{-\frac{\pi}{2}}^{\frac{\pi}{2}}\cos^{2(p'-1)}(\varphi)td\varphi\right)^{\frac{1}{p'}} dt \le c a^{-\frac{1}{p}}\|f\|_p.$$
If $p=1$,
$$|M_af(x)|\le \|f\|_1\frac{1}{A}\int_{B_a}\sqrt{t}\frac{1}{t}dt=\frac{c}{a}\|f\|_1.$$

In the parabolic case
$$|M_ag(x+h))-M_ag(x)|$$ $$\le \frac{1}{a\sqrt{\pi}}\int_{B_a} \left(\int_{|u|>R} (\cdot)+ \int_{|u|\le R}|g(x+h+2b\sqrt{t}u)-g(x+2b\sqrt{t}u)|e^{-u^2}du\right) dt.$$
Again in the first integral we refer to boundedness, in the second one to the uniform continuity of $g$. \\
Let $1 \le p<\infty$. Then
$$|M_af(x)|\le \frac{1}{\sqrt{\pi}} \|f\|_p\|e^{-u^2}\|_{p'}\frac{1}{A}\int_{B_a}(2b\sqrt{t})^{\frac{1}{p'}-1}dt \le c a^{-\frac{1}{2p}}\|f\|_p,$$
which is  \eqref{fy}.

\medskip

\begin{lemma} The translation defined in Example \ref{tm} is regular and $M_T=1$.
\end{lemma}

\proof
For boundedness see \cite[Theorem 2]{ms}.

Let our dense sets $\mathcal{E}=\mathcal{E}_i$, $i=1,2$ be the polynomials, that is
$$\mathcal{E}:=\left\{\sum_{k=0}^N a_k P_k^{(\lambda)}(\cos \vartheta): a_k\in \mathbb{R}, \ws \ws N\in\mathbb{N}\right\}.$$ Let $P\in\mathcal{E}$.
$$\|T^{r+h}P(\vartheta)-T^rp(\vartheta)\|_{\mu,p}=\|\sum_{k=0}^N a_k ((r+h)^k-r^k)P_k^{(\lambda)}(\cos \vartheta)\|_{\mu,p}$$ $$\le h c(N)\sum_{k=0}^N |a_k| \|P_k^{(\lambda)}(\cos \vartheta)\|_{\mu,p}=C(P) h,$$
which is (T1). (T2) can be shown similarly. Let $\nu(r)=\lambda(r)$, the Lebesgue measure on $[0,1)$ and $P\in\mathcal{D}$.
$$|M_{a,R}P(\vartheta+h)-M_{a,R}P(\vartheta)|=\left|\sum_{k=0}^N a_k \left(P_k^{(\lambda)}(\cos (\vartheta+h))-P_k^{(\lambda)}(\cos \vartheta)\right)\frac{1}{a}\int_0^a r^kdr\right|$$ $$=\left|-2h\lambda \sin \tilde{\vartheta}\sum_{k=0}^N \frac{a_k}{k+1}a^kP_{k-1}^{(\lambda+1)}(\cos \tilde{\vartheta})\right| \le c(P,a)|h|.$$
In view of \eqref{p} $P(r,\vartheta,\varphi)>0$, and if $r\le \frac{1}{8}$, say, then it is also bounded. Thus if $a\le \frac{1}{8}$, $1\le p<\infty$
$$|M_af(\vartheta)|\le\frac{1}{a}\int_0^a\|f\|_{\mu,p}\|P(r,\vartheta,\varphi)\|_{\mu,p'}dr\le c\|f\|_{\mu,p},$$
where $\|\cdot\|_{\mu,\infty}:=\|\cdot\|_{\infty}$.

\medskip

\medskip

\begin{lemma} The translation defined in Example \ref{tpl} is regular and bounded.
\end{lemma}

\proof
For boundedness see \cite[Theorem 2.2]{ch}.

The dense sets are the set of polynomials. By the notation \eqref{fi} it is as follows
$$\mathcal{E}_{1,2}=\mathcal{E}=\left\{\sum_{k=0}^n c_k \varphi_k^{(\alpha)}: \ws n\in\mathbb{N}, \ws c_k\in\mathbb{R}\right\}.$$
Let $f\in \mathcal{E}$, $f(x)=\sum_{k=0}^n c_k \varphi_k^{(\alpha)}(x)$. Then $T^tf(x)=\sum_{k=0}^n e^{-(2k+\alpha+1)t}c_k \varphi_k^{(\alpha)}(x)$. Thus
$$\left\|T^{t+h} f(x)-T^{t} f(x)\right\|_{p}=\left\|\sum_{k=0}^n e^{-(2k+\alpha+1)t}\left(e^{-(2k+\alpha+1)h}-1\right)c_k \varphi_k^{(\alpha)}(x)\right\|_p,$$
that is $h$ can be chosen appropriately depending on $f$.\\
Similarly, since $x\in B_R$
$$|M_{a,R}f(x+h)-M_{a,R}f(x)|=\left|\sum_{k=0}^n \frac{e^{-(2k+\alpha+1)a}-1}{(2k+\alpha+1)a}c_k \left(\varphi_k^{(\alpha)}(x+h)-\varphi_k^{(\alpha)}(x)\right)\right|$$
is small with an appropriate $h$ depending on $f$.\\
To prove the third property we consider the next estimation of the kernel function, see \cite[(2.4)]{fm}
$$0\ge W_t^\alpha(x,y)\le C\frac{e^{-c\frac{(x-y)^2}{t}}}{\sqrt{t}}.$$
Let $g\in L^p(I)$.
$$|M_{a,R}f(x)|\le \left|\int_0^{x-1}g(y)\frac{1}{a}\int_0^aW_t^\alpha(x,y)dtdy\right|+\left|\int_{x-1}^{x+1}(\cdot)\right|+\left|\int_{x+1}^\infty(\cdot)\right|=I+II+III.$$
$$II\le \|g\|_p\frac{c(p)}{\sqrt{a}}.$$
$$I,III\le \|g\|_p\frac{1}{\sqrt{a}}\left\|e^{-c\frac{(x-\cdot)^2}{a}}\right\|_{p'}.$$

\medskip

The consequence of Theorem \ref{T1} is the next one.

\begin{cor}
Let $K\subset L^p_\mu(I)$ be bounded. $K$ is precompact if and only if it is equivanishing and equicontinuous in mean with $T$ defined in Examples \ref{her}, \ref{th}, \ref{tp}, \ref{tm}, \ref{tpl} and \ref{jfv}.
\end{cor}

\section{Regularity with respect to an integral transform I - Pego type theorems}

In \cite{pe} Pego characterized the compact sets of $L^2(\mathbb{R}^d)$ by Fourier transform.  This result was extended to different spaces and to Abelian groups, see e.g. \cite{gk},  \cite{go} and the references therein. A Pego-type theorem by Laplace transform was proved in \cite{k}. These results based on standard translation. The extension of Pego's theorem to Bessel and Laguerre translations and transforms is given in \cite{hor}. The main ingredients of this type of theorems are a translation, a generated convolution and a corresponding integral transform.

\subsection{Convolution}
The study of convolution structures with general translations dates back to the 70s and it is a widely studied topic, see e.g. \cite{ag}, \cite{g}, \cite{fjk1}, \cite{fjk2}, \cite{ba} for Jacobi, \cite{gm} and \cite{ka}. These general convolutions, similarly to the standard ones are defined as follows. For $f$ and $g$ appropriate functions
$$f*g(t):=\int_IT^tf(x)g(x)d\mu(x).$$
All the convolutions in the listed papers are based on symmetric translations which are related to hyperbolic equations. This symmetry ensures algebraic structures with respect to the convolution in question since the symmetry of the kernel of the translation in its three variables implies the commutativity and associativity of the convolution, see the references above. In parabolic and elliptic cases the kernels of the translations are symmetric only in two variables thus the associativity fails.

To prove Pego-type theorems the next relations of convolution, integral transform and translation are necessary.
$$\mathcal{I}(f*g)=\mathcal{I}(f)\mathcal{I}(g),$$
and
$$\mathcal{I}(T^tf)(z)=\psi(t,z)\mathcal{I}(f)(z), \ws \ws T^t\mathcal{I}(f)(z)=\mathcal{I}(\psi f)(z),$$
where $\mathcal{I}$ is the integral transform and $\psi$ is an appropriate function. The translation is regular with respect to the integral transform if it possesses the properties above.

\subsection{Pego type theorems by standard convolution}

Pego type characterization of compactness can be given in cases of Examples \ref{th} and \ref{tp}.

In elliptic case we characterize compactness of sets in $L^2(\mathbb{R}_+)$ and consider cosine transformation, in parabolic case the sets are in $L^2(\mathbb{R})$ and the corresponding transformation is the standard Fourier one.\\
The transformation pairs are normalized as follows.
$$\mathcal{I}(f)(z):=\hat{f}(z)=\int_{\mathbb{R}} f(x)k(xz) dx, \ws \ws \ws \mathcal{J}(\hat{f})(x)=\frac{1}{2\pi}\int_{\mathbb{R}} \hat{f}(z)k(xz) dz,$$
where $\mathcal{I}=C$ or $\mathcal{I}=\mathcal{F}$, the cosine or the Fourier transform with $k(xz)=\cos(xz)$ or $k(xz)=e^{-ixz}$, respectively.
First the actions between translations and transformations are derived. Suppose that $f$ is in the Schwartz class and in elliptic case we assume that it is even.
$$T^tC(f)(z)=\frac{1}{\pi}\int_{\mathbb{R}}2\int_{\mathbb{R}_+} f(u)\cos u\eta du \frac{t}{(\eta-z)^2+t^2}d\eta$$ $$=2\int_{\mathbb{R}_+} f(u)\frac{1}{\pi}\int_{\mathbb{R}} \cos u\eta du \frac{t}{(\eta-z)^2+t^2}d\eta du=2\int_{\mathbb{R}_+} e^{-ut}f(u)\cos uz du$$ $$=C(e^{-|\cdot|t}f(\cdot))(z),$$
cf. \cite[1.2 (13)]{be}.
By similar calculation, supposing that $z\ge 0$,
$$C(T^tf(x))(z)=e^{-zt}C(f)(z).$$

In the parabolic case we have
$$T^t\mathcal{F}(f)(z)=\int_{\mathbb{R}}\int_{\mathbb{R}}\frac{e^{-\frac{(z-\eta)^2}{4b^2t}}}{2b\sqrt{\pi}\sqrt{t}}e^{-i\eta u}f(u)dud\eta=\int_{\mathbb{R}}T^t_z(e^{-i(\cdot)u})f(u)du$$ $$=\int_{\mathbb{R}}f(u)e^{-b^2tu^2}e^{-izu}du=\mathcal{F}\left(f(\cdot)e^{-b^2t(\cdot)^2}\right)(z).$$
On the other hand
\begin{equation}\label{pf}\mathcal{F}(T^tf(x))(z)=e^{-b^2tz^2}\mathcal{F}(f)(z).\end{equation}
As it is well-known,
$$\mathcal{I}: \ws L^p \longrightarrow L^{p'}, \ws \ws \ws 1\le p \le 2.$$
Let $K\subset L^p(I)$ ($1\le p \le 2$). Denote $\hat{K}\subset L^{'}p(0,\infty)$, $\hat{K}=\mathcal{I}(K)$. With this notation we have the next Pego-type theorem.

\medskip

\begin{theorem}\label{peg} Let $1\le p \le 2$, $K\subset L^p(I)$ a bounded set.\\ Let us consider Examples \ref{th} and \ref{tp}.\\ Assume that $K$ satisfies ${\bf P_a}$. Then $\hat{K}$ fulfils ${\bf P_{b_1}}$ and ${\bf P_{b_2}}$.\\
On the other hand if  $K$ satisfies ${\bf P_{b_1}}$, then $\hat{K}$ fulfils ${\bf P_a}$.
\end{theorem}

\proof Let $g(\cdot,t)=e^{-|\cdot|t}$ or $g(t,\cdot)=e^{-b^2t(\cdot)^2}$, and $f$ as above.
$$\|T^{t+h}\mathcal{I}(f)(z)-T^t\mathcal{I}(f)(z)\|_{p', x}= \|\mathcal{I}(g(\cdot,t+h)f(\cdot))(z)-\mathcal{I}(g(\cdot,t)f(\cdot))(z)\|_{p', z}$$ $$\le c \|g(x,t+h)f(x)-g(x,t)f(x)\|_p $$ $$\le \left(\int_{B_R}|f(x)|^p\left|g(x,t+h)-g(x,t)\right|^pdx\right)^{\frac{1}{p}}+2\left(\int_R^\infty |f(x)|^pdx\right)^{\frac{1}{p}}.$$
In view of ${\bf P_a}$, the second term is small. Since $x \in B_R$, $g(x,\cdot)$ is uniformly continuous. Thus considering that $K$ is bounded, ${\bf P_{b_1}}$ is satisfied for dense set uniformly and so for $\hat{K}$.

$$\left|M_{a,R}\mathcal{I}(f)(z+h)-M_{a,R}\mathcal{I}(f)(z)\right|$$ $$=\left|\frac{1}{A}\int_{B_a}\int_{\mathbb{R}}g(u,t)f(u)\left(k (u(z+h))-k(uz)\right)dud\nu(t)\right|$$ $$\le\int_{\mathbb{R}}\frac{1}{A}\int_{B_a}g(u,t)d\nu(t) |f(u)|\left|k (u(z+h))-k (uz)\right|du$$ $$\le c\left( \int _{B_R}h(u)|f(u)|\left|k(u(z+h))-k (uz)\right|du+\int _{\mathbb{R}\setminus B_R}(\cdot)\right)=c(III+IV),$$
where $h(u)=\frac{1}{A}\int_{B_a}g(u,t)d\nu(t)$ is uniformly bounded in $u$, and if $u\in \mathbb{R}\setminus B_R$,
$$h(u)\le \frac{c(a)}{|u|}$$
in both cases. Thus
$$IV\le \left(\int_{\mathbb{R}\setminus B_R} |f|^p\right)^{\frac{1}{p}}c(a)\left(\int_{\mathbb{R}\setminus B_R} u^{-p'}du\right)^{\frac{1}{p'}},$$
which is small by ${\bf P_a}$, considering that $p'>1$, and
$$III\le c\|f\|_p\|k(u(z+h))-k(uz)\|_{p',[0,R]}.$$
Since the interval is bounded, $h$ can be chosen such that the third factor be small and so, as above, $\hat{K}$ fulfils ${\bf P_{b_2}}$.

To prove the opposite direction we use the convolution below in cosine transform case, and the standard one in Fourier transform case. For cosine transform, supposing that $f$, $g$ are even
$$f*g(x):=\int_{\mathbb{R}}f(\xi)\frac{1}{2}(g(x+\xi)+g(x-\xi))d\xi.$$
Assume that $K$ fulfils ${\bf P_{b_1}}$.
Let $k_t(x)=K_{t,e}(x)$ or $k_t(x)=K_{t,p}(x)$, $x\in \mathbb{R}$, $y> 0$. Choose $R$ so large that $|\mathcal{I}(k_t)(\eta)|\le \frac{1}{2}$ if $\eta\in \mathbb{R}\setminus B_R$.
$$\left(\int_{\mathbb{R}\setminus B_R}|\mathcal{I}(f)(\eta)|^{p'}d\eta\right)^{\frac{1}{p'}}\le 2\left(\int_{\mathbb{R}\setminus B_R}|\mathcal{I}(f)(\eta)(1-\mathcal{I}(k_t)(\eta)|^{p'}d\eta\right)^{\frac{1}{p'}}$$ $$\le c \|f-f*k_t\|_p= c\|f -T^tf(x)\|_p.$$
Thus ${\bf P_{b_1}}$ implies ${\bf P_a}$.

\medskip

The next corollary implies that equicontinuity in Examples \ref{th} and \ref{tp} sense is equivalent with the standard one from this point of view.

\begin{cor} Let $K\subset L^2(I)$ be a bounded set and $\hat{K}$ the Fourier or the cosine transform of $K$, respectively. Then $K$ is precompact if and only if $K$ and $\hat{K}$ are equicontinuous in mean in sense of Examples \ref{tp} and \ref{th}, respectively.
\end{cor}

\medskip

\subsection{A Pego type theorem by general convolution}
In Laguerre and Bessel cases Pego-type theorems are proved in \cite{hor}. Below we prove a similar theorem with Jacobi transform. As the translation in the original and dual spaces are different, this situation is more complex than the Bessel case.

Jacobi transform is defined as it follows, see \cite[Proposition 3]{fj}. Recalling the measure $d\mu(x)=\frac{1}{\sqrt{2\pi}}w(x)dx$ (see \eqref{w} and \eqref{mu}), let $f\in L^2_\mu(\mathbb{R}_+)$ and $\lambda \in \mathbb{R}_+$.
\begin{equation}\mathcal{J}(f)(\lambda):=\hat{f}(\lambda)=\int_0^\infty f(x)\varphi_\lambda^{(\alpha,\beta)}(x)d\mu(x), \end{equation}
where the integral is convergent in $L^2_\nu(\mathbb{R}_+)$ with $d\nu(\lambda)=\frac{1}{\sqrt{2\pi}}v(\lambda)d\lambda$, where
\begin{equation} v(\lambda)=\left|\frac{2^{\varrho-i\lambda}\Gamma(i\lambda)\Gamma(\alpha+1)}{\Gamma\left(\frac{\varrho+i\lambda}{2}\right)\Gamma\left(\frac{\varrho+i\lambda}{2}-\beta\right)}\right|^{-2}.\end{equation}
$\mathcal{J}: f\mapsto \hat{f}$ is a linear and normpreserving map of $L^2_\mu(\mathbb{R}_+)$ onto  $L^2_\nu(\mathbb{R}_+)$. The inverse transform is given by
\begin{equation}\mathcal{J}^{-1}(f)(x)=f(x)=\int_0^\infty \hat{f}(\lambda)\varphi_\lambda^{(\alpha,\beta)}(x)d\nu(\lambda),
\end{equation}
where the integral is convergent in $L^2_\mu(\mathbb{R}_+)$.

Since $\mathcal{J}$ maps $L^1_\mu(\mathbb{R}_+)$ into $L^\infty(\mathbb{R}_+)$, see \cite[(3.2)]{fjk1}, Riesz-Thorin theorem ensures that $\mathcal{J}$ maps $L^p_\mu(\mathbb{R}_+)$ to  $L^{p'}_\nu(\mathbb{R}_+)$ continuously, where $1\le p\le 2$.

\medskip

As the range of the Jacobi transform is different from its domain, introduction of a dual translation is necessary, see \cite{fjk2}. It is given by the next formula.
\begin{equation}\label{Td}T_d^\eta g(\xi)=\int_0^\infty a(\xi,\eta,\zeta)g(\zeta)d\nu(\zeta),  \end{equation}
where the dual kernel is
$$a(\xi,\eta,\zeta)=\mathcal{J}(\varphi_\xi^{(\alpha,\beta)}\varphi_\eta^{(\alpha,\beta)})(\zeta)=\int_0^\infty \varphi_\xi^{(\alpha,\beta)}(x)\varphi_\eta^{(\alpha,\beta)}(x)\varphi_\zeta^{(\alpha,\beta)}(x)d\mu(x),$$
see \cite[(4.14)]{fjk2}. That is the kernel is symmetric in its three variables. Moreover
$$a(\xi,\eta,\zeta)\ge 0, \ws (\xi,\eta,\zeta)\in\mathbb{R}^3; \ws\ws \int_0^\infty a(\xi,\eta,\zeta)d\nu(\zeta)=1,$$
see \cite[Theorem 4.4 and (4.17)]{fjk2}. This ensures, in standard way, that
\begin{equation}\label{tdf}\|T_d^\eta g\|_{p,\nu}\le \|g\|_{p,\nu}, \ws \ws 1\le p \le \infty.\end{equation}

\medskip

Subsequently we need the next properties (see also \cite[Section 2]{git}).

\begin{lemma} Let $f\in C_0(\mathbb{R}_+)$. Then
\begin{equation}\label{JT}\mathcal{J}(T^t f)(\lambda)=\varphi_\lambda^{(\alpha,\beta)}(t)\mathcal{J}(f)(\lambda)\end{equation}
and
\begin{equation}\label{TJ}T_d^\eta\hat{f}(\lambda)=\mathcal{J}(f\varphi_\eta^{(\alpha,\beta)})(\lambda).\end{equation}
\end{lemma}

\proof By \cite[(4.2)]{fjk1} $$\int_0^\infty K(x,t,z)\varphi_\lambda^{(\alpha,\beta)}(z)d\mu(z)=\varphi_\lambda^{(\alpha,\beta)}(t)\varphi_\lambda^{(\alpha,\beta)}(x).$$ Thus
$$\mathcal{J}(T^t f)(\lambda)=\int_0^\infty \int_0^\infty f(z)K(x,t,z)d\mu(z)\varphi_\lambda^{(\alpha,\beta)}(x)d\mu(x)$$ $$=\int_0^\infty f(z)\varphi_\lambda^{(\alpha,\beta)}(t)\varphi_\lambda^{(\alpha,\beta)}(z)d\mu(z),$$
where by the assumption on $f$ Fubini theorem can be applied and the last integral is convergent. It proves \eqref{JT}.

Similarly, in view of \cite[(4.16)]{fjk2} $$\int_0^\infty a(\xi,\eta,\zeta)\varphi_\zeta^{(\alpha,\beta)}(x)d\nu(\zeta)=\varphi_\xi^{(\alpha,\beta)}(x)\varphi_\eta^{(\alpha,\beta)}(x),$$ which implies that
$$T_d^\eta\hat{f}(\lambda)=\int_0^\infty \int_0^\infty f(x)a(\lambda,\eta,\zeta)\varphi_\zeta^{(\alpha,\beta)}(x)d\mu(x)d\nu(\zeta)$$ $$=\int_0^\infty f(x)\varphi_\eta^{(\alpha,\beta)}(x)\varphi_\lambda^{(\alpha,\beta)}(x)d\mu(x)$$
which, with the remark above, ensures \eqref{TJ}.

\medskip

Before stating the next theorem we introduce the corresponding convolution. For appropriate functions
$$f*g(t):=\int_0^\infty T^tf(x)g(x)d\mu(x)$$
and
$$\mathcal{J}(f*g)=\mathcal{J}(f)\mathcal{J}(g),$$
furthermore the convolution is commutative and associative, see \cite[(5.3), (5.4)]{fjk1}.

\begin{theorem}\label{pegg} Let $1\le p \le 2$, $K\subset L^p_\mu(\mathbb{R}_+)$ a bounded set. Let $\hat{K}:=\mathcal{J}(K)$.\\  Let us consider Example \ref{jfv}.\\
If $K$ satisfies ${\bf P_a}$, then $\hat{K}$ fulfils ${\bf P_{b_1}}$ and ${\bf P_{b_2}}$.\\
If  $K$ satisfies ${\bf P_{b_1}}$, then $\hat{K}$ fulfils ${\bf P_a}$.
\end{theorem}

\proof Let $f\in K$. In view of \eqref{TJ}
$$\|T_d^{\eta+h}\hat{f}(\lambda)-T_d^\eta\hat{f}(\lambda)\|_{p',\nu}=\|\mathcal{J}(f(\varphi_{\eta+h}^{(\alpha,\beta)}-\varphi_\eta^{(\alpha,\beta)}))\|_{p',\nu}\le c\|f(\varphi_{\eta+h}^{(\alpha,\beta)}-\varphi_\eta^{(\alpha,\beta)})\|_{p,\mu}.$$
According to \cite[Theorem 2, (ib)]{fj} if $\eta \in \mathbb{R}$,
$$\frac{\partial}{\partial \eta}\varphi_{\eta}^{(\alpha,\beta)}(x)\le K(1+x)^2e^{-\varrho x}.$$
Thus, since  $\varrho>0$, $\|T_d^{\eta+h}\hat{f}(\lambda)-T_d^\eta\hat{f}(\lambda)\|_{p',\nu}\le c(\varrho)h\|f\|_{p,\mu}$, so ${\bf P_{b_1}}$ is satisfied for $\hat{K}$.

Repeating the standard arguments, for an arbitrary $\delta>0$ take a function $g_\delta$ such that $\mathrm{supp}g_\delta \subset [-\delta,\delta]$, $\int_{\mathbb{R}_+}g_\delta(t)d\mu(t)=1$, $\hat{g}_\delta>0$, cf. \cite[Lemma 4.2]{fjk2}. As $\hat{g}_\delta(\lambda)$ tends to zero when $\lambda$ tends to infinity, choose $R$ so large that $\hat{g}_\delta(\lambda)<\frac{1}{2}$ if $\lambda>R$. Then
$$\left(\int_R^\infty |\hat{f}|^{p'}d\nu\right)^{\frac{1}{p'}}\le \left(\int_R^\infty |\hat{f}(1-\hat{g}_\delta)|^{p'}d\nu\right)^{\frac{1}{p'}}
\le c(p)\|f-f*g_\delta\|_{p,\mu}$$ $$=c(p)
\left(\int_0^\infty \left|\int_0^\infty (f(x)-T^x f(t))g_{\delta}(t)d\mu(t)\right|^pd\mu(x)\right)^{\frac{1}{p}} $$ $$\le c(p)\left(\int_0^\infty \int_0^\infty \left|f(x)-T^tf(x)\right|^p d\mu(x)g_\delta(t)d\mu(t)\right)^{\frac{1}{p}}$$ $$\le c(p)\sup_{0\le|t|\le\delta} \|f(x)-T^tf(x)\|_{p,\mu}.$$

\section{Modulus of smoothness and K-functional - Regularity with respect to an integral transform II}

As it is well-known, the modulus of smoothness generated by the standard translation is equivalent with the Peetre's K-functional, see e.g. \cite[page 171]{dl}. This property is extended to Bessel translation, see \cite{p}. Below we derive the same equivalence for heat-diffusion semigroup, where the semigroup property implies arguments very similar to the standard case, and for Example \ref{tp}, where the corresponding integral transform proves to be the right tool.\\
Recalling the notation of \eqref{des} we can define the Sobolev space generated by $D_{(2)}$ as follows.
\begin{equation}W_{D,r,\mu}^p=\{f(x)\in L^p_\mu :  D_{(2),x}^l f\in L^p_\mu, \ws l=1, \dots , r\}.\end{equation}
Then the corresponding K-functional is
$$K_r(f,t,L^p_\mu, W_{D,r,\mu}^p):=K_r(f,t)_p=\inf_{g\in W_{D,r,\mu}^p}\{\|f-g\|_{p,\mu}+t\|(D_{(2)}^r g\|_{p,\mu}\}.$$
The moduli of smoothness generated by general translations are defined as follows.

\begin{defi} Let
$$\Delta_tf(x):=T^tf(x)-f(x), \ws \ws \ws \Delta_t^rf(x):=\Delta_t(\Delta_t^{r-1}f)(x).$$
Let $1\le p\le \infty$. The $p$-modulus of smoothness is
$$\omega_r(f,t)_p:=\sup_{0<h\le t}\|\Delta_h^rf\|_{p, \mu}.$$
\end{defi}

\medskip

\remark (R4)
\begin{equation}\label{D}\Delta_t^rf(x)=(I-T^t)^rf=\sum_{k=0}^r(-1)^{r-k}\binom{r}{k}(T^t)^kf(x),\end{equation}
where
$(T^t)^k=T^t\circ \dots \circ T^t$. For the standard and the heat-diffusion translation semigroups $(T^t)^k$ becomes $T^{kt}$.

\medskip

\begin{theorem}\label{mck} There are positive constants $M=M(r,p)$ independent of $f$ such that for Examples \ref{her} and \ref{tp}
$$\frac{1}{M}\omega_r(f,t)_p \le  K_r(f,t^{cr})_p\le M\omega_r(f,t)_p,$$
$r\in\mathbb{N}$, $1\le p\le \infty$, $c=\mathrm{ord}(D_{(1)})$.
\end{theorem}

\medskip

\remark

\noindent (R5) In  our parabolic cases $c=1$, in Bessel case $c=2$, cf. \cite{p}.

\medskip

\proof  By standard arguments
$$\omega_r(f,t)_p \le \omega_r(f-g,t)_p+\omega_r(g,t)_p \le 2^r\|f-g\|_p + t^r\|D_{(2)} g\|_p.$$
In view of \eqref{delta} and Remark (R4),
$$\|\Delta_t^rg\|_p\le t^r\|D_{(2)} g\|_p.$$

The second estimation needs different methods in the two cases.

First we deal with the heat semigroup. Let us introduce the next notation.
$$g_{r,t}(x):=f(x)+(-1)^{r+1}r^r\int_{(0, \frac{1}{r})^r}\Delta_{t(\sum_{i=1}^rz_i)}^rf(x)dz_1\dots dz_r$$ \begin{equation}\label{gth}=r^r\int_{(0, \frac{1}{r})^r}\sum_{k=1}^r(-1)^{k+1}\binom{r}{k}u(x,kt(\sum_{i=1}^rz_i))dz_1\dots dz_r.\end{equation}
Thus, according to \eqref{hde}
\begin{equation}\label{Lgh} D_{(2),x}^rg_{r,t}(x)=\frac{r^r}{t^r}\sum_{k=1}^r(-1)^{k+1}\binom{r}{k}\frac{1}{k^r}\Delta_{\frac{kt}{r}}^rf(x).\end{equation}
Considering the operator norm, \eqref{gth} and \eqref{Lgh} $g_{r,t} \in W^p_{r,L}$. By Minkowski's inequality
$$t^r\|D_{(2),x}^rg_{r,t}(x)\|_p\le (2r)^r\omega_r(f,t)_p.$$
Recalling the definition of $g_{r,t}$, as above
$$\|f-g_{r,t}\|_p= \left\|r^r\int_{(0, \frac{1}{r})^r}\Delta_{t(\sum_{i=1}^ru_i)}^rf(x)du_1\dots du_r\right\|_p\le M(r,p)\omega_r(f,t)_p.$$

To prove the second inequality for Example \ref{tp} we apply the integral transform method again.
According to \eqref{D} and \eqref{pf}
$$(\mathcal{F}(\Delta_t^rf))(z)=(1-e^{-b^2tz^2})^r\hat{f}(z).$$
Let $I_1:=I\cap [-1,1]$ and $I_2:=I\setminus (-2,2)$. Let us define
$$\eta(x)= \left\{\begin{array}{ll} 1, \ws  x\in I_1,\\0, \ws x\in I_2\\ \in C_0^\infty \mbox{otherwise},\end{array}\right.$$
such that $0\le \eta \le 1$. Let $\varepsilon>0$ and define $\eta_\varepsilon(x):=\eta\left(\frac{x}{\varepsilon}\right)$. Let
$$\varrho_\varepsilon:=\mathcal{F}^{-1}\left(\eta_\varepsilon\right).$$
It is enough to take an $f\in L^1\cap L^p$. Define
$$g:=g_\varepsilon:=(\varrho_\varepsilon*f).$$
Let Let $0<t<\frac{1}{b^2\varepsilon^2}$. Then we have
$$\mathcal{F}(D_{(2)}^rg_\varepsilon)(z)=(-1)^rb^{2r}z^{2r}\eta_\varepsilon(z)\hat{g}(z)$$ $$=\frac{1}{t^r}\frac{t^r(-1)^rb^{2r}z^{2r}\eta_\varepsilon(z)}{(1-e^{-b^2tz^2})^r}(\mathcal{F}(\Delta_t^rg))(z)=\frac{1}{t^r}\frac{t^r(-1)^rb^{2r}z^{2r}\eta(bz\sqrt{t})}{(1-e^{-b^2tz^2})^r}(\mathcal{F}(\Delta_t^rg))(z).$$
Let
$$h(u)= \left(\frac{-u^2}{1-e^{-u^2}}\right)^r\eta(u).$$
Then
$$D_{(2)}^rg_\varepsilon=\frac{1}{t^r}\Delta_t^rg\mathcal{F}^{-1}(h(bz\sqrt{t})).$$
Since
$$\|\mathcal{F}^{-1}(h(bz\sqrt{t}))\|_1=\|\mathcal{F}^{-1}(h)\|_1=:c,$$
we have
\begin{equation}\label{2tag}\|D_{(2)}^rg_\varepsilon\|_p\le c\frac{1}{t^r}\|\Delta_t^rg\|_p.\end{equation}
We proceed similarly to estimate the first term of $K_r(f,t^{r})_p$.
$$\mathcal{F}(f-g_\varepsilon)(z)=\left(1-\eta_\varepsilon\right)\hat{g}(z)=\frac{1-\eta_\varepsilon}{(1-e^{-b^2tz^2})^r}\widehat{\Delta_{t}^rg}(z).$$
Let
$$\frac{1}{(1-x)^r}=\sum_{k=0}^\infty a_kx^k=Q_N(x)+\sum_{k=N+1}^\infty a_kx^k, \ws \ws |x|<1.$$
We decompose
$$\frac{1}{(1-e^{-b^2tz^2})^r}=\left(\frac{1}{(1-e^{-b^2tz^2})^r}-Q_N\left(e^{-b^2tz^2}\right)\right)+Q_N\left(e^{-b^2tz^2}\right),$$
that is
$$\mathcal{F}(f-g_\varepsilon)(z)=U_1+U_2$$ $$=\left(1-\eta_\varepsilon\right)Q_N\left(e^{-b^2tz^2}\right)\widehat{\Delta_{t}^rg}(z)$$ $$+\left(1-\eta_\varepsilon\right)\left(\frac{1}{(1-e^{-b^2tz^2})^r}-Q_N\left(e^{-b^2tz^2}\right)\right)\widehat{\Delta_{t}^rg}(z).$$
Let $\varepsilon=\frac{1}{b\sqrt{t}}$ and
$$k_N(z)=(1-\eta(z))\left(\frac{1}{1-e^{-z^2}}-Q_N(z)\right).$$
Then
$$\mathcal{F}^{-1}U_2=\mathcal{F}^{-1}\left(k_N(bz\sqrt{t})\right)*\Delta_{t}^rg.$$
If $N$ is large enough we have
\begin{equation}\label{els}\|\mathcal{F}^{-1}U_2\|_p\le \|\mathcal{F}^{-1}k_N\|_1\|\Delta_{t}^rg\|_p.\end{equation}
Turning our attention to the first term we have
$$\mathcal{F}^{-1}U_1=\sum_{k=0}^Na_k\left(\left(T^t\right)^k(\Delta_{t}^rg)-\varrho_{\frac{1}{bz\sqrt{t}}}*(\left(T^t\right)^k(\Delta_{t}^rg)\right).$$
Since $T^t$ is bounded
\begin{equation}\label{mas}\|\mathcal{I}^{-1}U_1\|_p\le C(Q)(1+\|\varrho _1\|_1)\|\Delta_{t}^rf\|_p.\end{equation}
Thus according to \eqref{els} and \eqref{mas}
$$\|f-g_\varepsilon\|_p\le c \omega_r(f,t)_p.$$
Recalling \eqref{2tag} the second inequality is proved.

\medskip

\remark (R6)

This integral transform method works because considering the two expressions $\mathcal{I}(\Delta_tf)=(1-\psi(t,\cdot))\hat{f}$ and $\mathcal{I}D_{(2)}f=h(\cdot)\hat{f}$ we have
$$h(\cdot) \sim 1-\psi(t,\cdot)$$
around zero. This property is the second regularity property of translation with respect to the integral transform. This is the situation e.g. in Bessel case, see \cite{p}. For instance Example \ref{th} does not possess second regularity property with respect to the cosine transform.

\medskip

\noindent \small{Department of Analysis, Institute of Mathematics,\newline
Budapest University of Technology and Economics \newline
 M\H uegyetem rkp. 3., H-1111 Budapest, Hungary.

\medskip

 g.horvath.agota@renyi.hu}
\end{document}